\documentclass{article}
\usepackage{amsthm}  
\usepackage{amssymb}
\usepackage{amsmath}
\usepackage{mathptmx}

\newcommand{\mbb}[1]{\mathbf{#1}}

\newcommand{\mrm}[1]{\mathrm{#1}}
\newcommand{\C}{\mbb{C}}

\newcommand{\N}{\mbb{N}}

\newcommand{\Pro}{\mbb{P}}

\newcommand{\Aut}{\mrm{Aut~}}
\newcommand{\GAut}{\mrm{GAut}}
\newcommand{\fn}[3]{{#1:#2\rightarrow #3}}

\theoremstyle{definition}
\newtheorem{definition}{Definition}[section]
\newtheorem{remark}[definition]{Remark}

\newtheorem{example}[definition]{Example}
\theoremstyle{plain}
\newtheorem{theorem}[definition]{Theorem}

\newtheorem{lemma}[definition]{Lemma}
\newtheorem{proposition}[definition]{Proposition}

\date{}		

\begin{document}

\title{Maximally symmetric stable curves\footnote{2000MSC: 14H37}}
\author{Michael A. van Opstall\and R\u{a}zvan Veliche}

\maketitle

\begin{abstract}
We prove a sharp bound for the automorphism group of a stable curve of a 
given genus and describe all curves attaining that bound.
\end{abstract}

All curves are defined over $\C$ and projective.

A well-known result of Hurwitz \cite{h} states that the maximal order of the 
automorphism group of a smooth curve of genus $g$ is $42(2g-2)$. This
bound is not attained in every genus; for example, the maximal order
of the automorphism group of a smooth genus two curve is forty-eight, attained 
by the curve with affine equation $y^2=x^5+x$. In genus three, the bound of
168 is achieved by the famous {\em Klein quartic} $y^7=x^3-x^2$
 (given in homogenous 
coordinates by
$x^3y+y^3z+z^3x$). A curve attaining the Hurwitz bound is known as a 
{\em Hurwitz curve}; a great deal is known about these curves and the 
corresponding automorphism groups.

Exercise 2.26 of \cite{hm:mc} asks if this bound holds for {\em stable 
curves}. A stable curve is a curve with nodal singularities and finite
automorphism group. Gluing three copies of the elliptic curve with 
$j$-invariant zero (which has six automorphisms which fix a point) to a
copy of the projective line yields a genus three curve with $6^3\cdot6$
automorphisms, breaking the Hurwitz bound. In genus two, gluing the 
aforementioned elliptic curve to itself yields a curve with seventy-two
automorphisms - not breaking the Hurwitz bound - but more symmetric than any
smooth genus two curve.

The goal of this work is to give a sharp bound for the automorphism group 
of a stable curve of genus $g$ and more-or-less describe all curves 
attaining this
bound. Since the moduli space of stable curves is locally the quotient of
a smooth germ by the automorphism group of a stable curve, this bound gives
some measure of how singular this moduli space can be near the boundary. 

\section{Geometric preparation}

Denote by $E$ the most symmetric elliptic curve (that with $j$-invariant 
zero) in what follows. We will also use $\Pro^1$ somewhat abusively to
denote a smooth rational curve. $\Pro^1$ will be coordinatized as the
Riemann sphere. Recall that to a stable curve one may associate a {\em dual
graph} which is a (weighted) graph (possibly with multiple edges and loops)
with one vertex for each component of the curve 
(labelled with its genus) and vertices are connected by edges if the 
corresponding components meet at a node. Self-intersecting curves lead to
loops in this graph. Given a vertex $v$, we will call the number of edges
connecting $v$ to a vertex corresponding to a $\Pro^1$ the {\em rational
valence} of $v$. {\em Elliptic valence} is defined similarly.

Automorphisms of stable curves come from two sources: automorphisms of
their components which preserve or permute the nodes properly, and certain
automorphisms of the dual graph. Not every graph automorphism is induced by
an automorphism of the curve. For example, $n$ points on $\Pro^1$ can be
permuted at most dihedrally. Once another point is required to be fixed,
$n$ additional points may be permuted at most cyclically. After fixing two
points (say zero and infinity), $n$ points may still be permuted cyclically 
(the $n$th roots of unity). All attaching of curves to copies of $\Pro^1$
will be tacitly done in the most efficient way: placing several isomorphic
branches of the curve at roots of unity, and making other attachments at
zero and infinity as not to disrupt the cyclic symmetry. 

\begin{definition}
An automorphism of a dual graph $G$ will be called {\em geometric} if it is
induced by an automorphism of the corresponding stable curve. The group
of such automorphisms will be denoted $\GAut~G$.
\end{definition}

We will call the vertices of the dual graph corresponding to the rational
curves {\em rational vertices} and those corresponding to elliptic curves
{\em elliptic vertices}. Finally, a vertex in a graph meeting a single edge 
will always be called a {\em leaf}, whether or not the graph in question is a 
tree. As usual, the corresponding components of the curve are called 
{\em tails}.

\begin{definition}
A {\em maximally symmetric stable curve} is a stable curve whose automorphism
group has maximal order among all stable curves of the same genus.
\end{definition}

Note that this definition makes sense: a stable curve of genus $g$ has at 
most $2g-2$ components, each with normalization of genus at most $g$ (note
the vulgarity of this bound). Therefore there is a bound for the automorphism
group of a genus $g$ curve of $[42(2g-2)]^{2g-2}(2g-2)!$.

In this section we shrink the class of curves which need to be considered. An
initial lemma will be essential in what follows:

\begin{lemma}\label{stab}
The dual graph of the stabilization of a nodal curve of genus $g>1$ with
no tails has at least as many automorphisms as the dual graph of the original 
curve.
\end{lemma}

\begin{proof}
First note that if a vertex is deleted in stabilizing a graph, all
vertices in the orbit of this vertex are also deleted.

We proceed by induction on the number of vertices in the graph. There are no
unstable curves of genus two or higher with a graph with a single node, so 
there is nothing to prove in this case. 

The only possibility is that there is a vertex $v$ of valence two corresponding
to a rational curve. Proceed in both directions along geodesics away from
$v$ until vertices $u$ and $w$ are reached which are stable (i.e., not 
deleted in stabilizing). Such vertices exist since the curve has genus at 
least two,
and they may coincide. Whenever $v$ is moved by an automorphism, the arc
from $u$ to $w$ is moved with it, and obviously conversely. Therefore, 
replacing the entire arc from $u$ to $w$ with a single edge from $u$ to $w$ 
(including the case of replacing the arc with a loop from $u=w$ to itself)
does not decrease the automorphism group of the curve.
\end{proof}

\begin{lemma}\label{smcmp}
A maximally symmetric curve has only smooth components.
\end{lemma}

\begin{proof}
Let $C$ be any stable curve, and suppose that $C_1$ is a component with nodes.
Replacing $C_1$ in $C$ with its normalization drops the genus of $C$ by
the number of nodes of $C_1$. For each pair of points of the normalization
lying over a node, choose one and glue a copy of $E$ to it. This makes up
the genus deficit. The automorphisms of each copy of $E$ multiply the order
of the automorphism group by six. The automorphisms of $C_1$ correspond to
those of its normalization which permute the nodes appropriately. Therefore,
the normalization of $C_1$ has more automorphisms than $C_1$, and fewer of 
these are
killed off by gluing elliptic curves to only one of each pair of points
over a node than identifying these points. If the normalization is genus 
zero, and $C_1$ has only one node, the normalized component will be collapsed 
to keep the curve stable, and $E$ reattached at
the point of attachment of the rational curve.

Finally, if multiple isomorphic copies of $C_1$ occur in $C$ resulting in
a symmetry of the dual graph which induces automorphisms of $C$, replace
{\em each} copy of $C_1$ by this construction to maintain the symmetry. We call
this {\em maintaining graph symmetry}, and omit explicit mention of it in
the future.
\end{proof}

\begin{figure}[ht]
\begin{center}
$\begin{array}{c@{\hspace{.5in}}c@{\hspace{.5in}}c}
\begin{picture}(40,40)
\put(11,30){\vector(1,0){8}}

\put(5,20){\circle*{3}}
\put(5,20){\line(0,1){10}}
\put(5,30){\circle*{3}}
\put(5,33){\circle{6}}

\put(25,20){\circle*{3}}
\put(25,20){\line(0,1){10}}
\put(25,30){\circle*{3}}
\put(25,30){\line(0,1){10}}
\put(25,40){\circle*{3}}
\put(30,37){\makebox{\scriptsize 1}}
\end{picture}
&
\begin{picture}(40,60)
\put(11,30){\vector(1,0){8}}

\put(5,20){\circle*{3}}
\put(5,20){\line(0,1){10}}
\put(5,30){\circle*{3}}
\put(8,28){\makebox{\scriptsize 3}}

\put(25,20){\circle*{3}}
\put(25,20){\line(0,1){10}}
\put(25,30){\circle*{3}}
\put(28,25){\makebox{\scriptsize 0}}
\put(25,30){\line(0,1){10}}
\put(25,40){\circle*{3}}
\put(28,37){\makebox{\scriptsize 1}}
\put(25,30){\line(-1,1){10}}
\put(15,40){\circle*{3}}
\put(18,37){\makebox{\scriptsize 1}}
\put(25,30){\line(1,1){10}}
\put(35,40){\circle*{3}}
\put(38,37){\makebox{\scriptsize 1}}
\end{picture} 

&

\begin{picture}(50,45)
\put(15,30){\vector(1,0){8}}

\put(10,20){\circle*{3}}
\put(10,20){\line(0,1){10}}
\put(10,30){\circle*{3}}
\put(12,28){\makebox{\scriptsize 1}}
\put(10,30){\line(0,1){10}}
\put(10,40){\circle*{3}}

\put(30,20){\circle*{3}}
\put(30,20){\line(0,1){10}}
\put(30,30){\circle*{3}}
\put(25,28){\makebox{\scriptsize 0}}
\put(30,30){\line(0,1){10}}
\put(30,40){\circle*{3}}
\put(30,30){\line(1,0){10}}
\put(40,30){\circle*{3}}
\put(38,24){\makebox{\scriptsize 1}}

\end{picture}
\\
\mbox{\bf Lemma \ref{smcmp}} & \mbox{\bf Lemma \ref{ratell}}
& \mbox{\bf Lemma \ref{nomult}}
\end{array}$
\end{center}
\caption{Illustrating Lemmas \ref{smcmp}-\ref{nomult}}
\label{fig-smre}
\end{figure}

\begin{lemma}\label{ratell}
There exists a maximally symmetric stable curve whose components are all
$\Pro^1$ or $E$.
\end{lemma}

\begin{proof}
By Lemma \ref{smcmp}, all components may be taken smooth, and it is clear that
replacing an elliptic curve less symmetric than $E$ with $E$ only helps.

If $C_1$ is a genus $h>1$ component of $C$, replace it with a rational curve,
attach a rational tail to this, and then glue $h$ copies of $E$ to this second
rational curve. Ignoring graph symmetry,
$C_1$ contributed at most $42(2h-2)$ automorphisms to $C$, whereas the new
construction contributes at least $2\cdot 6^h$ (the two comes from the fact
that there are at least two copies of $E$ which may be permuted, since $h>1$).

There is a technical issue: $2\cdot 6^2$ is not bigger than $84$, but we have
already seen that there is no smooth genus two curve with more than 
$2\cdot 6^2$
automorphisms. Also, in the case that $C$ itself is a smooth genus two curve,
this construction leads to a non-stable curve. This is fine - the maximally
symmetric genus two curve is two copies of $E$ glued together, which 
satisfies the conclusion of the lemma.
\end{proof}

\begin{lemma}\label{nomult}
A maximally symmetric stable curve's components are all
copies of $\Pro^1$ or $E$; the dual graph of such a curve has no multiple 
edges, and its leaves are elliptic, and other vertices are rational.
\end{lemma}

\begin{proof}
Apply the constructions of Lemmas \ref{smcmp} and \ref{ratell}. Suppose there
is a copy of $E$ which is not a tail of the curve. Then $E$ is attached in
at least two points, so replacing $E$ with a $\Pro^1$ and gluing $E$ to this
$\Pro^1$ does not decrease the number of automorphisms, and makes $E$ a
tail. 

By Lemma \ref{smcmp} there are no loops in the graph. Suppose two vertices
are connected by $n$ edges. Since we may assume at this point that all
elliptic components are tails, these vertices must be rational. Replace
the multiple edge by a vertex joined by two edges to the former endpoints
of the multiple edge. Add a rational tail to this new vertex, and arrange
$n-1$ copies of $E$ as tails around this rational vertex. The curve may need
to be stabilized, but we have seen that in this case stabilization will
not affect automorphisms.

If $n=3$ and the entire curve
is two $\Pro^1$ attached in three points, this
construction does not result in a stable curve, but again, this is an 
exceptional case, and we know that the maximally symmetric genus two curve
satisfies the conclusions of the lemma.

The configuration of two rational curves connected by $n$ nodes and connected
some other way to the rest of the curve contributes (excluding graph symmetry)
at most $2n$ automorphisms. This construction replaces this with a 
configuration contributing $6^{n-1}$ automorphisms.
\end{proof}

\begin{figure}[ht]
\begin{center}
$\begin{array}{c@{\hspace{.5in}}c}
\begin{picture}(32,40)
\put(11,30){\vector(1,0){8}}

\put(7,20){\circle*{3}}
\put(7,20){\line(0,1){10}}
\put(7,30){\circle*{4}}
\put(0,27){\makebox{\scriptsize 0}}
\put(6,30){\line(0,1){10}}
\put(8,30){\line(0,1){10}}
\put(7,40){\circle*{4}}
\put(0,37){\makebox{\scriptsize 0}}

\put(27,20){\circle*{3}}
\put(27,20){\line(0,1){10}}
\put(27,30){\circle*{3}}
\put(31,27){\makebox{\scriptsize 0}}
\put(27,30){\line(0,1){10}}
\put(27,40){\circle*{3}}
\put(31,37){\makebox{\scriptsize 1}}
\end{picture}

&

\begin{picture}(80,60)
\put(36,30){\vector(1,0){8}}

\put(10,20){\circle*{3}}
\put(5,18){\makebox{\scriptsize 0}}
\put(10,20){\line(1,0){20}}
\put(30,20){\circle*{3}}
\put(32,18){\makebox{\scriptsize 0}}
\put(10,20){\line(1,2){10}}
\put(30,20){\line(-1,2){10}}
\put(20,40){\circle*{3}}
\put(24,38){\makebox{\scriptsize 0}}

\put(55,30){\circle*{3}}
\put(55,30){\line(0,1){10}}
\put(55,40){\circle*{3}}
\put(58,37){\makebox{\scriptsize 0}}

\put(55,30){\line(-1,-1){10}}
\put(45,20){\circle*{3}}
\put(47,18){\makebox{\scriptsize 0}}

\put(55,30){\line(1,-1){10}}
\put(65,20){\circle*{3}}
\put(67,18){\makebox{\scriptsize 0}}

\put(55,30){\line(1,0){10}}
\put(65,30){\circle*{3}}
\put(68,27){\makebox{\scriptsize 1}}
\end{picture} 

\\
\mbox{\bf Lemma \ref{nomult}} & \mbox{\bf Lemma \ref{tree}}
\end{array}$
\end{center}
\caption{Illustrating Lemmas \ref{nomult}-\ref{tree}}
\label{fig-tree}
\end{figure}

\section{Breaking cycles}

This section contains the main part of the reduction: we may assume that
the dual graph of a maximally symmetric stable curve is a tree. To achieve
this, we need to break cycles in the graph, producing a new graph with
more automorphisms. In one case, this is easy. In this entire section we
assume that all of the reductions from the previous section have been carried
out: we are now studying simple graphs whose interior vertices all correspond
to smooth rational curves and whose leaves are copies of $E$.

\begin{definition}
A cycle in a graph is called {\em isolated} if it shares no edge with any
cycle in its orbit under the action of the automorphism group of the graph.
\end{definition}

\begin{lemma}\label{tree}
The dual graph of a maximally symmetric stable curve $C$ has no isolated
edge-transitive cycles.
\end{lemma}

\begin{proof}
Such a cycle of $n$ rational
curves contributes one to the genus and contributes at most $2n$ automorphisms
(dihedral symmetry).
Replace this cycle of $n$ curves with a ``wheel'' whose hub is a rational 
curve with $n$ rational ``spokes'' connected at the $n$th roots of unity and a 
copy of $E$ attached at zero. This possibly drops graph symmetry by a factor of
two (reflections in the dihedral group are not included in this case, because
our automorphisms must fix the point of attachment of the spoke), 
but multiplies automorphisms by six due to the introduction of a copy of
$E$.
\end{proof}

\begin{remark}
The assumption that the cycle is isolated is necessary so that the construction
can be carried out on every cycle in the orbit.
\end{remark}

The following proposition is obvious, but we state it for ease of reference:

\begin{proposition}
\label{vertexclasses}
There are at most two orbits of the automorphism group among the vertices of 
an edge-transitive graph.
\end{proposition}

\begin{definition}
Suppose $G$ is a graph and $G'$ is a subgraph. Then $\GAut_{G'}G$ will
denote the group of geometric automorphisms of $G$ which fix $G'$.
\end{definition}

In what follows, a graph will be called {\em optimal} if its geometric 
automorphism group is maximal among dual graphs of stable curves of a
given genus.

\begin{proposition}(Edge transitive graphs are not optimal)
\label{estimate.edge.transitive}
Let $C$ be a stable curve with dual graph $G$.
If $G$ is an edge-transitive graph with valence at least three at each of its 
vertices, there exists a curve $C'$ whose dual graph is a tree with elliptic 
leaves such that 
$|\Aut C'|>|\Aut C|$.
\end{proposition}

\begin{proof}
By \ref{vertexclasses}, there are two orbits of vertices, $O(v)$ and $O(w)$, 
with valences
$e_v$ and $e_w$, respectively; 
set $n_v=|O(v)|$ and $n_w=|O(w)|$. 
Then $n_ve_v=n_we_w$ and the genus of $G$ is 
$g=\frac{n_v}{2}(e_v-2)+\frac{n_w}{2}(e_w-2)+1$. The number of vertices of $G$ 
is $n=n_v+n_w$. The case in which there exists a single orbit $O(v)$ is dealt 
with exactly as in the case where $e_v=e_w$.

We will bound $|\GAut~G|$ using a sequence of trees that ``grow'' and 
eventually include all the vertices of $G$ (spanning trees). 

Start with a vertex $v$; denote by $T_0$ the tree consisting of $v$ alone. 
There are $n_v$ choices for $T_0$.

For $T_1$ take $T_0$ and all vertices (which are in the orbit of $w$) 
adjacent with $T_0$. Then $|\GAut_{T_0}T_1|\leq 2e_v$ (the automorphisms of 
$T_1$ fixing $T_0$ at most dihedrally permute the $e_v$ edges around 
$v$).

Assume the trees $T_0,\dots, T_i$ have been constructed. If $T_i$ spans $G$, 
we stop. If not, there is a vertex of $T_i$, call it $x$, such that 
at least one of its neighbors is not in $T_i$. Let $T_{i+1}$ be the span of
$T_i$ and $x$. Set $n_i$ to be the number of vertices in the tree $T_i$,
and $e$ the valence of $x$ (which equals either $e_v$ or $e_w$).

We have the following possibilities:
\begin{enumerate}
\item if $e\geq 4$ and $x$ has at least three neighbors in $T_i$, or if 
$e=3$ and $x$ has two neighbors in $T_i$, then $|\GAut_{T_i}T_{i+1}|=1$. 
In both cases $n_{i+1}\geq n_i+1$.
\item if $e\geq 4$ and $x$ has at most two neighbors in $T_i$, or if $e=3$ 
and $x$ has exactly one neighbor in $T_i$, then 
$|\GAut_{T_i}T_{i+1}|\leq 2$.
In this case $n_{i+1}\geq n_i+(e-2)$ (respectively $n_{i+1}=n_i+2$).
\end{enumerate}

This process will terminate after a finite number of steps, since $G$ has 
finitely many vertices.

Denoting by $s_v$ the number of times the second possibility occurs with 
$x\in O(v)$ and by $s_w$ the number of times the second possibility occurs 
with $x\in O(w)$, we have $n_w\geq e_v+s_v(e_v-2)$ (respectively 
$n_w\geq 3+2s_v$ when $e_v=3$) and $n_v\geq 1+s_w(e_w-2)$ (respectively 
$n_v\geq 1+2s_w$ when $e_w=3$).  At the same time, it is clear that 
$|\GAut~G|\leq n_v\cdot 2e_v\cdot 2^{s_v+s_w}$.

A curve of genus $g$ whose dual graph is a tree with elliptic leaves will 
have at least $6^g$ automorphisms.
We want to show that $6^g> 3\cdot 2n_ve_v\cdot 2^{s_v+s_w}$. This means  
$6^{\frac{n_v}{2}(e_v-2)+\frac{n_w}{2}(e_w-2)}> n_ve_v 2^{s_v+s_w}$. 

We have several cases to consider, depending on the valencies $e_v$ and $e_w$.

\begin{enumerate} 
\item $e_v=e_w=3$. Then $s_v\leq \frac{n_w-3}{2}$ and 
$s_w\leq \frac{n_v-1}{2}$, so it is sufficient to prove 
$6^{\frac{n_v}{2}+\frac{n_w}{2}}> 3n_v 2^{\frac{n_w-3}{2}+\frac{n_v-1}{2}}$ 
or $\sqrt{6}^n> \frac{3}{4}n\sqrt{2}^n$, which is clearly true.
\item $e_v=3<e_w$. Then $s_w\leq \frac{n_v-1}{e_w-2}$ and 
$s_v\leq \frac{n_w-3}{2}$; the inequality to prove becomes
$6^{\frac{n_v}{2}+\frac{n_w}{2}(e_w-2)}\geq 
3n_v 2^{\frac{n_w-3}{2}+\frac{n_v-1}{e_w-2}}$. 
Now $n_v=n_w\frac{e_w}{3}$, so the inequality becomes (removing the index $w$): 
$6^{\frac{ne}{6}+\frac{n}{2}(e-2)}\geq 
ne 2^{\frac{n-3}{2}+\frac{ne-3}{3(e-2)}}$ or yet 
$6^{\frac{2ne}{3}-n}\geq 
ne\cdot 2^{\frac{5n}{6}-\frac{3}{2}+\frac{2n-3}{3(e-2)}}$ 
(since $\frac{ne-3}{3(e-2)}=\frac{n}{3}+\frac{2n-3}{3(e-2)}$). 
This is true, since $e>3$, by the inequality 
$6^{\frac{2ne}{3}-n}\geq ne\cdot 2^{2n}$, 
or $(\frac{3}{2})^n\cdot 6^\frac{2n(e-3)}{3}\geq ne$. For $n\geq 2$, 
$(\frac{3}{2})^n\geq n$, and $6^\frac{4(e-3)}{3}\geq 6^{e-3}> e$ (since $e>3$).
For $n=1$, $\frac{3}{2}\cdot 6^{\frac{2}{3}(e-3)}> \frac{3}{2}3^{e-3}>e$ since 
$e>3$. So in this case we are done also.
\item $e_v, e_w>3$. Then $s_w\leq \frac{n_v-1}{e_w-2}$ and 
$s_v\leq \frac{n_w-e_v}{e_v-2}$. The inequality to prove becomes 
$6^{\frac{n_v}{2}(e_v-2)+\frac{n_w}{2}(e_w-2)}\geq 
n_ve_v\cdot 2^{\frac{n_v-1}{e_w-2}+\frac{n_w-e_v}{e_v-2}}$; 
since $e_v,e_w>3$, this is implied by 
$6^{n_ve_v-n_v-n_w} > n_ve_v 2^\frac{n_v+n_w}{2}$. 
Using $n_w=\frac{n_ve_v}{n_w}$ and dropping the index $v$ we get 
$6^{ne}>ne\cdot (6\sqrt{2})^{n+\frac{ne}{e_w}}$. Now $6\sqrt{2}<9$ and 
$e_w\geq 4$, so it is enough to show that 
$6^{ne}>ne\cdot 3^{2n+\frac{ne}{2}}$, 
and yet again $(2\sqrt{3})^{ne}>ne\cdot 3^{2n}$. Since $2n\leq \frac{ne}{2}$,
this is implied by $(2\sqrt{3})^{ne}>ne\cdot (\sqrt{3})^{ne}$ which becomes 
$2^{ne}>ne$ which is finally clear.
\end{enumerate}

\end{proof}


\begin{proposition}(Collapsing edge-transitive subgraphs)
\label{finalstep}
Assume that $G$ is a dual graph which consists of  
an edge-transitive graph $H$ with rational vertices $v_1,\ldots,v_n$ 
with a tree $T_i$ attached at each vertex $v_i$ which is either
degenerate (i.e. consists of $v_i$ only) or has only elliptic leaves and 
otherwise rational vertices. To avoid trivialities, assume that $H$ is
not a tree. Then $G$ is not an optimal graph.
\end{proposition}

\begin{proof}
The point is that $H$ has at least one cycle, and this forces less than 
desirable symmetry in $G$.

By \ref{vertexclasses} we know that there are at most two orbits of vertices 
in $H$. Assume that there are exactly two, $O(v)$ and $O(w)$; the other case 
is similar (practically identical proof using $e_v=e_w$).

Denote by $n_v$ and $n_w$ the order of $O(v)$ and $O(w)$ in $H$, and by 
$e_v$ and $e_w$ their respective valence in $H$. Then $n_ve_v=n_we_w$.

If $e_v=e_w=2$, $H$ is an isolated edge-transitive cycle in $G$, and thus 
$G$ cannot be optimal.

The genus of $H$, $g(H)=\frac{n_v}{2}(e_v-2)+\frac{n_w}{2}(e_w-2)+1\geq 2$. 

If all the $T_i$ are degenerate, $G$ is simply $H$ and using 
\ref{estimate.edge.transitive} we know that there exists an optimal tree with 
strictly more automorphisms than $G$; so we are done in this case.

If some $T_i$ are not degenerate, the genus of $H$ is smaller than the genus 
of $G$, so by induction there exists an optimal tree $T$ with 
$|\GAut~T|>|\GAut~H|$ (note that \ref{estimate.edge.transitive} applies
to graphs with some vertices of valence two in light of \ref{stab}). Detach 
the non-degenerate isomorphic $T_i$ (including 
the edge that connects their root to the vertex $v_i$) from the vertices of 
$H$, pair them two-by-two around a new root, in the end connecting these 
roots of pairs to a new root $V$; connecting $V$ to the root of $T$ by an 
edge leaves the overall genus unchanged, and yields at least a two-fold 
increase in automorphisms (actually, a $2^\frac{n-1}{2}$, if $n$ is the number 
of the $T_i$ detached); this increase is due to the extra freedom allowed on 
the $T_i$, which can be swapped independently of the vertices $v_i$. The
reattaching can 
at most decrease the symmetry of $T$ by a factor of two, so as soon as one 
orbit of 
isomorphic trees has at least two elements, we get the desired strict increase 
in symmetry; in particular, $G$ was not optimal.

If there exists a unique tree $T_i$ which is non-degenerate, necessarily 
there exists an orbit of vertices of $H$ with a single vertex in it. Since 
$H$ was edge-transitive, the condition that $H$ is not a tree implies that 
$H$ is ``tree-like'' (of diameter two), but with multiple edges. 
From the breaking of multiple edges earlier in the proof, we get at least a 
threefold increase in symmetry by shooting out elliptic tails in place of 
multiple edges; again in particular $G$ was not optimal.
\end{proof}

So now all we need to do is show that starting with a graph which is not a 
tree, with elliptic leaves and otherwise rational vertices, there is a graph 
of the same genus, satisfying the properties of the previous lemma, with at 
least as many automorphisms. 
Then the previous lemma shows that the original graph could not 
have been optimal.

\begin{lemma}(Pre-valence reduction)
\label{prep}
We may assume that an optimal graph of genus $g$ has the following property: 
around each vertex there are at most three orbits of edges, and at most one 
orbit of two or more edges.
\end{lemma}

\begin{proof} 
Assume that in an optimal graph we have a vertex $v$, necessarily rational 
(by previous reductions), around which there exist either four or more orbits, 
or at least two orbits of edges, each with at least two edges (ending at $v$) 
in it. 

Let us establish some notation: arrange the orbits of edges around $v$ in 
decreasing order of their sizes; so we have 
$O_1,\ldots,O_k,O_{k+1},\ldots,O_{k+l}$, where 
$|O_1|\geq |O_2|\geq\cdots |O_k|\geq 2>1=|O_{k+1}|=\cdots=|O_{k+l}|$ ($k$ 
or $l$ may be zero).

So we must show that there exists an optimal graph in which $k\leq 1$ and 
$k+l\leq 3$. 

Note that if $k+l\geq 4$, there are no automorphisms permuting the edges 
around $v$, and if $k+l=3$ there is at most one (non-trivial) orbit being 
cyclically (half-dihedrally) permuted by the automorphisms fixing $v$.

We will perform the following operations on the given graph:

\begin{enumerate}
\item detach all edges from around $v$, keeping track of their orbits; 
\item replace $v$ by a path of rational vertices 
$v_1-v_2-\ldots-v_k-v_{k+1}-\ldots-v_{k+l-1}$;
\item attach the orbit $O_i$ to $v_i$ for $1\leq i\leq k+l-2$, and attach the 
orbits $O_{k+l-1}$ and $O_{k+l}$ to $v_{k+l-1}$. 
\end{enumerate}

The operation above should be done simultaneously at all vertices in the 
orbit of $v$, so as not to lose the initial graph symmetry; the fact that 
these vertices are in the same orbit implies that the same partition of edges 
is repeated around each such vertex, and thus the same insertion of the new 
rational path can be done everywhere. Note that the genus of the graph has not 
been changed, as one vertex has been replaced by $k+l-1$ vertices and 
$k+l-2$ edges.

Note that there is a map from the new graph to the old graph. Prescribing that 
the newly introduced paths will go (orientation-preserving) only onto another 
similar (newly introduced) path lifts distinct automorphisms of the original 
graph to distinct automorphisms of the new graph. Thus we have preserved or
increased the order of the automorphism group (i.e. the new graph is also 
optimal), and the claim of the lemma is established.
\end{proof}


We will need the following lemma in changing graphs into those of the type in 
\ref{finalstep}.

\begin{lemma}(Bound on graphs with fixed vertices)
\label{bound.fix.vertices}
Let $G$ be a simple graph of genus $g\geq 1$ with elliptic leaves and 
otherwise 
rational vertices; assume that the valence at each rational vertex is at least 
three, except possibly at some rational vertices $v_1,\dots,v_k$ where it may 
be 
two. Let $H$ be a subgroup of $\GAut~G$ which fixes the vertices $v_1,\ldots,
v_k$ and acts at most cyclically (i.e. not fully dihedrally) on the edges 
around them; then $2|H| < |\GAut~T|$ for some $T$ a tree with $g$ elliptic 
tails.
\end{lemma}

\begin{proof}
The lemma is obviously true for $g=2$, inspecting case-by-case (and noting 
that the vertices of valence two do not bring any extra symmetry).

Assume that $H$ fixes only one rational vertex $v$. The orbit of this vertex 
by $H$ is trivial. By valence reasons, removing this vertex and the 
$k\geq 2$ edges around it (which can only be in an orbit by themselves) will 
yield a graph $G'$ of genus $g-k+1$, with valence at least two at any vertex, 
but will not decrease the automorphism group other than by at most a factor of 
$k$ (due to the action on the removed edges being at most cyclic); in other 
words, the subgroup of $H$ fixing the edges around $v$ has index at most $k$ 
in $H$. Moreover, the automorphisms fixing the edges around $v$ will fix their 
opposite ends in $G'$. Now either $G'$ has genus one (with our 
assumptions on valencies on $G$, that means a cycle, which would be fixed by 
any automorphism fixing the edges around $v$), or has genus at least two. In 
either case one can replace inductively $G'$ by a tree $T'$ with at least 
twice as many automorphisms (in the case of the cycle, one may simply use 
another elliptic tail), arrange another $k-1$ elliptic tails around a root, 
link this root to the root of $T'$ at $\infty$ and fix $0$ as well. At most a 
dihedral symmetry is lost at the root of the tree replacing $G'$ in this 
manner, but $6^{k-1}$ is gained through the elliptic tails. Overall one gains 
at least a factor of $\frac{6^{k-1}}{k}>2$ in the automorphism group, which 
is what is needed.

Now if $H$ fixes several rational vertices $v_1,\dots,v_k$, 
we may use the previous case to bound the larger subgroup fixing only one of 
the $v_i$ and get the desired bound.

\end{proof}

\begin{theorem}
The optimal graph of genus $g$ is a tree. More precisely, for any graph $G$ 
of  genus $g$ which is not a tree there exists a tree $T$ (with $g$ elliptic 
tails) such that $|\GAut~T|>|\GAut~G|$.
\end{theorem}

\begin{proof}
Begin with an optimal graph which is pre-valence reduced (see \ref{prep}). 
Recall that we already 
know that an optimal tree should not have isolated cycles.

Choose an edge $e$ such that:
\begin{itemize}
\item the order of its orbit is smallest among all orbits of edges,
\item in case this order is at least two, require that the orbit of one of 
its ends be the smallest among 
the possible choices for $e$.
\end{itemize}

The goal of this choice is to control the valence of the graph left by 
removing $e$ and the edges in its orbit.

To fix notations for the remainder of the proof, denote by $v$ the end of $e$ 
with the smallest orbit, and 
by $w$ the other end; denote by $n_v$ the order of the orbit of $v$, with a 
similar notation for $w$.


\begin{lemma}(Structure of a minimal edge orbit)
\label{min.edge.orbit} 
Assume $e$ has a unique representative around $v$. Then one and only one of 
the following cases can occur:
\begin{enumerate}
\item $e$ has a unique representative around $w$ as well.
\item All the edges around $w$ are in $O(e)$.
\end{enumerate}
\end{lemma}

\begin{proof} 
Assume that $e$ would have at least another representative around $w$, and 
that there would exist another edge $f\notin O(e)$; 
then $v$ and $w$ are not 
in the same orbit, and $|O(e)|\geq 2n_w$; in the same time, $|O(f)|\leq n_w$ 
(by pre-valence reduction)
contradicting the choice of $e$.
\end{proof}

We will prove the Theorem inductively. Let us establish some more notation: 
the connected components of 
$G'=G\setminus O(e)$ will be denoted by $G_1,\dots,G_k$; the graph obtained 
from $G$ by contracting the components $G_1,\dots, G_k$ to vertices will be 
called 
$G''$ (this could have multiple edges, i.e. be non-simple); some optimal tree 
with elliptic tails and the same genus as $G_i$ ($G''$) will be called $T_i$ 
(respectively $T''$).

The basic idea is to look at the connected components of $G'$, replace them 
inductively (if necessary) by optimal trees, then reconnect the trees at 
their roots to form a common tree. Some care needs to be taken with this 
procedure because some components may be degenerate (isolated vertices), and 
some symmetry might be lost in the individual trees (those with full dihedral 
symmetry around their root) when they are connected (by an edge ending at 
their root) to some other trees. However, the Lemma \ref{bound.fix.vertices} 
shows that the last part is not a real concern.

\begin{proposition}
\label{valence}
With the given choice of $e$:
\begin{enumerate}
\item If around a vertex $v$ there are three edges in $O(e)$, then all edges
ending at $v$ are in $O(e)$.
\item Removing the edges in the orbit of $e$ from the graph $G$ cannot leave a 
vertex with valence one. 
\end{enumerate}
\end{proposition}

\begin{proof}
For the first part, if $f$ is an edge ending at $v$ not in $O(e)$,
$|O(e)|\geq \frac32n_v$, while $|O(f)|\leq n_v$, contradicting
the choice of $e$.

Now, if $O(e)$ has only one element  ending at $v$ and $w$, then we are 
done: by stability, there must be at least two other edges, not in the orbit 
of $e$ (unless one or both of these is a tail, in which case the valence left 
is zero), around both $v$ and $w$.

So assume that there are at least two edges in the orbit of $e$ ending at $v$.
If there are at least three such edges, then the first part of this
Proposition says that all the edges around $v$ are in $O(e)$, so $v$ would have
valence zero in $G'$. If there are precisely two edges in $O(e)$ around $v$,
then stability of $G$ and Lemma \ref{min.edge.orbit} show that $e$ is not
unique around $w$ as well.

\begin{itemize}
\item If $O(e)$ has only two representatives around $w$, then the existence of 
an edge-transitive isolated cycle formed by edges in $O(e)$ 
with vertices in $O(v)$ and $O(w)$ is immediate. 
\item If $O(e)$ has at least three representatives around $w$, then 
actually all edges around $w$ must be in $O(e)$, by the first part. 
Then, denoting by $f$ (one of) the extra edge(s) around $v$, we have 
that $|O(e)|\geq \frac{2n_v+3n_w}{2}>n_v\geq |O(f)|$ which would lead 
to a contradiction.
\end{itemize}
\end{proof}


Note that for $G''$ to be a tree, the only possibility is that $e$ is unique 
in its orbit around $v$, and that all $G_i$'s containing a vertex in $O(w)$ 
contain a unique such vertex, and no vertex in $O(v)$.

There are two possibilities for $G'$: it is connected or disconnected.

{\bf Case 1}: $G'$ is connected. Then there are no isolated vertices left in 
$G'$, and all vertices have valence at least two; using \ref{valence}, $e$ 
must be unique in its orbit around both $v$ and $w$. The automorphism group 
of $G'$ has order at least that of $G$ (examining the movement of vertices). 
The genus of $G'$ is at least one less than that of $G$; stabilizing $G'$ 
does not decrease the number of automorphisms and preserves the genus, except 
when $G'$ was a cycle -- but then it would be isolated in its orbit, and not
adjacent to any edge-transitive cycles, so $G$ would not be optimal by an 
argument similar to \ref{tree}.
  
By induction, some optimal graph in a genus (at least) one less is a tree 
$T'$ and $|\GAut~T'|>|\GAut~G'|$; now compensate for the loss of genus by 
attaching the necessary number of elliptic tails to the root of $T'$. If 
$T'$ had dihedral symmetry at the root, the loss of it (factor of two) is 
easily compensated since each elliptic tail gives a factor of six increase in 
automorphisms. But then we would reach a contradiction to the fact that $G$ 
was optimal.

{\bf Case 2}: $G'$ is disconnected. 

If all the components $G_i$ are single vertices (i.e. when all the edges 
around both 
$v$ and $w$ are in $O(e)$), we are done because one of the following holds:
\begin{itemize}
\item All these vertices are rational, in which case $G$ is an 
edge-transitive graph with only rational vertices (so there are at most two 
orbits of vertices in it); the Lemma \ref{estimate.edge.transitive} shows 
that these are not optimal, i.e. this case cannot occur with our choice of $e$.
\item All of these vertices are elliptic, in which case $G$ was the dual 
graph of two elliptic curves meeting in a node, thus a tree.
\item Some of the vertices are elliptic and some rational. In this case it is 
clear that there can only be one rational vertex and all the elliptic vertices 
were connected by $e$ and its translates to it; thus $G$ was already a tree 
(of diameter two).
\end{itemize}

So we may assume that some component $G_i$ is not an isolated point and thus 
must be of positive genus; then the genus of $G''$ is strictly less than the 
genus of $G'$ and we are in the situation described in \ref{min.edge.orbit} 
(i.e. $O(e)$ does not exhaust all edges around both ends of $e$).

Note that we have the following formula bounding $|\GAut~G|$: 
$|\GAut~G|\leq |\GAut^GG''|\cdot\prod\limits_{i=1}^k |\GAut_{v,w}G_i|$, where 
$\GAut_{v,w}G_i$ is the group of automorphisms of $G_i$ fixing the ends of
edges in $O(e)$, and $\GAut^GG''$ is the group of automorphisms of $G''$ 
induced by
those of $G$. This follows because the automorphisms fixing $O(e)$ 
automatically fix vertices in each component $G_i$.

If $G$ has no cycles, it is a tree and we are done. We will assume implicitly 
from now on that $G$ has a cycle and show that it is not optimal.

Construct a graph $H$ in the following way: take $G''$ and attach to its 
vertices the trees $T_i$ via a new edge, at their roots. Note that the genus 
of $H$ is the same as the genus of $G$, which is
$\sum\limits_{i=1}^k g(G_i)+g(G'')$. 

When does this construction provide a stable graph with more symmetry than 
the original $G$? Rather: how much is the symmetry of the graph affected by 
this construction?

Lemma \ref{bound.fix.vertices} shows that in replacing the $G_i$ with 
$T_i$, if necessary, no symmetry is lost. It is similarly clear that fixing 
only the root (as opposed to any other vertex) of a tree yields the maximum 
number of automorphisms.

Let $f$ denote the root of the tree $T_i$ (which may be an edge, see section 
4).
Note that $|\GAut~H|\geq |\GAut_GG''|\cdot\prod\limits_{i=1}^k|\GAut_fT_i|$ 
actually. Thus $H$ must be optimal, too. Now \ref{finalstep} finishes the 
proof of the fact that $G$ was not optimal.
\end{proof}

\section{Valence reduction and statement of the Main Theorem}

\begin{definition}
If the rational valence at a vertex is $r$ and the elliptic valence is $e$,
then we say the {\em valence} of this vertex is $(r,e)$. This will cause no
confusion with the usual use of the word valence.
\end{definition}

\begin{lemma}\label{valred}
In an optimal graph, the valence of each rational vertex may only be one of: 
(0,3), (0,4), (0,5), (3,0), (4,0), (5,0), 
(1,2), (1,3), (2,1) or (3,1).
\end{lemma}

\begin{proof}
``Smaller'' valences are ruled out by stability of the curve. Suppose there
is a point of valence $(n,0)$ with $n\geq 6$. If all of the branches from
this vertex are mutually non-isomorphic, then the vertex can be replaced
by a chain of rational vertices and various branches distributed in a way
that decreases the valence into the allowed range. This will not affect
automorphisms. On the other extreme, if all vertices are isomorphic, if $n=2k$,
we can replace the vertex with a chain of $k$ rational vertices and attach
the branches to this pairwise. This replaces dihedral symmetry of order
$4k$ with $k$ involutions, plus a global involution of the chain - at least
$2^{k+1}$ automorphisms. If $k=2$ this does not affect the order of the
automorphism group, but if $k\geq 3$, the order increases. If $n$ is odd, 
a similar procedure grouping three branches together and otherwise pairing
branches works similarly. In the intermediate cases where there are some 
isomorphic branches but not all branches are isomorphic, split the vertex
into a chain of rational vertices, one for each isomorphism class of branches,
reattach the branches, and apply the arguments above. 

The cases $(0,n)$ are handled similarly.

For a vertex of the form $(1,n)$ with $n>3$, pair the elliptic leaves as
much as possible and connect the resulting two leaf branches (and possibly
a single leaf) to the vertex. This will transfer excess elliptic valence to
rational valence, which can then be distributed as above. In this case, since
there is one branch which cannot possibly be isomorphic to the others, to 
make symmetry maximal, this one branch (the rational one) should be glued to
the origin of $\Pro^1$, and the elliptic branches as symmetrically as possible
at roots of unity. But an automorphism which fixes the origin can only permute
the roots of unity cyclically, so we can reduce the valence further than
in the previous cases.

The cases $(n,1)$ are handled similarly.

Finally, vertices of valence $(r,e)$ with $r$ and $e$ both larger than allowed
here are dealt with by adding two branches to the vertex in question, 
distributing the rational branches of the original vertex around one, and
the elliptic vertices around another. This reduces to two vertices of
types $(1,e)$ and $(r+1,0)$, which have already been dealt with.
\end{proof}

Note that we cannot do better than this Lemma in general: if we have a 
vertex of valence $(5,0)$, and all the branches are isomorphic, the 
contribution to symmetry near this vertex is dihedral of order ten. On the
other hand, if we split this into two vertices of valence $(4,0)$ and $(3,0)$
connected by an edge, the connecting edge corresponds to the components
of the curve being glued together; the most symmetric option is to glue
the two curves at their origins, and the branches at roots of unity. But 
gluing the origins together makes the symmetry of the roots of unity only
cyclic, so we have six automorphisms only. Splitting into two branches
with two leaves and a single leaf gives at most eight automorphisms.

\begin{definition}
A stable curve is {\em simple} if 
\begin{enumerate}
\item its components are all copies of $E$ or $\Pro^1$,
\item its dual graph is a tree with all leaves elliptic and
all other vertices $\Pro^1$;
\item the valence of each vertex is no greater than five, and the elliptic
valence is no greater than three;
\end{enumerate}
\end{definition}

The previous reductions show that there exists a maximally symmetric 
stable curve of genus $g$ is simple. The ``geometric'' contribution to the 
automorphism
group of a genus $g$ simple curve is $6^g$, and the rest comes from certain
automorphisms of the graph.

Under the assumption of simplicity, elliptic components are
distinguished from rational components by occuring as leaves on the tree,
so we need not count automorphisms of the dual graph as a weighted graph. 
Therefore, the problem of finding a maximally symmetric stable curve has been
reduced to finding a maximally geometrically symmetric graph of a certain
type. 

\begin{figure}[ht]
\begin{center}
$\begin{array}{c@{\hspace{0.5in}}c@{\hspace{0.5in}}c}
\begin{picture}(20,48)
\put(5,24){\circle*{3}}
\put(5,24){\line(1,0){10}}
\put(15,24){\circle*{3}}
\end{picture}

&

\begin{picture}(48,48)
\put(24,24){\circle*{3}}

\put(24,24){\line(1,0){8}}
\put(24,24){\line(-1,0){8}}
\put(24,24){\line(0,1){8}}
\put(24,24){\line(0,-1){8}}

\put(24,8){\circle{16}}
\put(8,24){\circle{16}}
\put(24,40){\circle{16}}
\put(40,24){\circle{16}}

\put(20,5){\makebox{$2^n$}}
\put(4,21){\makebox{$2^n$}}
\put(20,37){\makebox{$2^n$}}
\put(37,22){\makebox{$a$}}
\end{picture}

&

\begin{picture}(72,48)
\put(36,24){\circle*{3}}
\put(8,24){\circle*{3}}
\put(64,24){\circle*{3}}

\put(36,24){\line(1,0){28}}
\put(36,24){\line(-1,0){28}}
\put(36,24){\line(0,1){8}}
\put(8,24){\line(0,1){8}}
\put(8,24){\line(0,-1){8}}
\put(64,24){\line(0,1){8}}
\put(64,24){\line(0,-1){8}}

\put(8,8){\circle{16}}
\put(8,40){\circle{16}}
\put(36,40){\circle{16}}
\put(64,8){\circle{16}}
\put(64,40){\circle{16}}

\put(4,5){\makebox{$2^n$}}
\put(4,37){\makebox{$2^n$}}
\put(60,5){\makebox{$2^n$}}
\put(60,37){\makebox{$2^n$}}
\put(33,38){\makebox{$b$}}
\end{picture}

\end{array}$
\end{center}
\caption{Illustrating three of the cases of the Main Theorem}
\label{mtfig}
\end{figure}
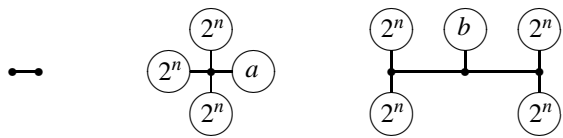

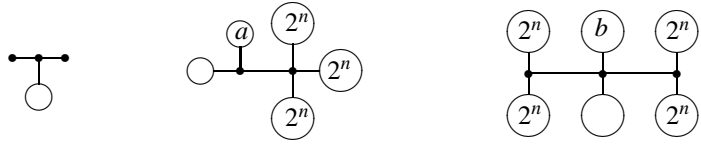
\begin{figure}[ht]
\begin{center}
$\begin{array}{c@{\hspace{0.5in}}c@{\hspace{0.5in}}c}

\begin{picture}(40,45)
\put(20,15){\circle{10}}
\put(20,20){\line(0,1){10}}
\put(20,30){\circle*{3}}
\put(20,30){\line(-1,0){10}}
\put(20,30){\line(1,0){10}}
\put(10,30){\circle*{3}}
\put(30,30){\circle*{3}}
\end{picture}

&

\begin{picture}(85,56)
\put(5,25){\circle{10}}
\put(10,25){\line(1,0){10}}
\put(20,25){\circle*{3}}
\put(20,25){\line(0,1){9}}
\put(20,25){\line(1,0){20}}
\put(20,39){\circle{10}}
\put(18,37){\makebox{$a$}}
\put(40,25){\circle*{3}}
\put(40,25){\line(0,1){10}}
\put(40,25){\line(1,0){10}}
\put(40,25){\line(0,-1){10}}
\put(58,25){\circle{16}}
\put(53,22){\makebox{$2^n$}}
\put(40,43){\circle{16}}
\put(37,40){\makebox{$2^n$}}
\put(40,7){\circle{16}}
\put(37,4){\makebox{$2^n$}}
\end{picture}

&

\begin{picture}(72,48)
\put(36,24){\circle*{3}}
\put(36,24){\line(0,-1){8}}
\put(36,8){\circle{16}}

\put(8,24){\circle*{3}}
\put(64,24){\circle*{3}}

\put(36,24){\line(1,0){28}}
\put(36,24){\line(-1,0){28}}
\put(36,24){\line(0,1){8}}
\put(8,24){\line(0,1){8}}
\put(8,24){\line(0,-1){8}}
\put(64,24){\line(0,1){8}}
\put(64,24){\line(0,-1){8}}

\put(8,8){\circle{16}}
\put(8,40){\circle{16}}
\put(36,40){\circle{16}}
\put(64,8){\circle{16}}
\put(64,40){\circle{16}}

\put(4,5){\makebox{$2^n$}}
\put(4,37){\makebox{$2^n$}}
\put(60,5){\makebox{$2^n$}}
\put(60,37){\makebox{$2^n$}}
\put(33,38){\makebox{$b$}}
\end{picture}

\end{array}$
\end{center}
\caption{Attachment rules for attaching appendages of types 1-3, from left to right, to a tree (denoted by an empty circle).}
\label{attach}
\end{figure}

\begin{theorem}[Main Theorem]\label{main}
In genus $g$, a maximally symmetric curve may be constructed as follows:
\begin{enumerate}
\item if $g=2$, two copies of $E$ glued at a point;
\item if $g=3\cdot 2^n+a$ for some $a<2^n$ and $n\geq 0$, the graph is three 
binary trees
attached to a central node which also has a branch which is a maximally
symmetric graph of a genus $a$ curve;
\item if $g=4\cdot 2^n+b$ for some $b<2^{n+1}$ (but $b\neq 2^n$) and 
$n\geq 1$, the graph is four binary
trees, two attached on one side of a chain of length three (two if $b=0$)
and two on the other, with a maximally symmetric graph of a genus $b$ curve
attached to the middle vertex of the chain.
\item if $g=5\cdot 2^n$, the graph is five binary trees attached to a common
vertex.
\end{enumerate}
(See Figure \ref{mtfig}).
Here when $a=1$ or $b=1$, a maximally symmetric genus one curve is assumed
to be a copy of $E$. See Figure \ref{attach} for how the appendages $a$
and $b$ are attached (stabilization may be required after attaching certain
appendages).
See \ref{form} for a computation of the orders of the automorphism groups in
each of these four cases.
\end{theorem}

\begin{remark}
In Section 5, it will be convenient to think of the third case given in
the theorem as a single binary tree with $2^{n+2}$ nodes attached to an
appendage $b$, but the picture of four smaller binary trees is a little 
clearer.
\end{remark}

Note that the fourth case is almost a subcase of the third, except the 
appendage balloon is a binary tree, so more automorphisms occur from collapsing
the graph.

The proof of this Theorem is graph-theoretic and deserves its own section.

\section{Proof of the Main Theorem}

A tree has either an edge or vertex which is invariant under the action
of the automorphism group (see Corollary 2.2.10 of \cite{s:tr}). If a
vertex is invariant, it will be called a {\em root} of the tree. If an
edge is invariant, call it a {\em virtual root}. If no confusion is possible,
either one will be called a {\em root}. We will actually need something
a little stronger:

\begin{lemma}
If $G$ is a tree of finite diameter $n$, then all geodesics of length $n$ have
the same middle vertex if $n$ is odd, or a common middle edge if $n$ is even.
\end{lemma}

\begin{proof}
This is Exercise 2.2.3 in \cite{s:tr}.
\end{proof}

Such a vertex (edge) will be called an {\em absolute (virtual) root}. Since
$G$ is a tree, given any other vertex $v$ of $G$, there is a unique edge 
adjacent to $v$ along the geodesic to an absolute (virtual) root. Let $G_v$ 
denote the subtree of $G$ obtained as follows: remove the edge adjacent
to $v$ along the geodesic to an absolute root; $G_v$ is the connected 
component containing $v$ of what remains.

In the rest of this section, $G$ is always assumed to be the dual graph of a 
simple stable curve.

\begin{definition}
Given a graph $G$, let $V(G)$ be the set of vertices of $G$ and $E(G)$ be the
set of edges. 
Define $\fn{o_V}{V(G)}{\N}$ by $o_V(v)=|O(v)|$. When we speak of an 
automorphism
acting on an edge, the edge is assumed to be oriented (that is,
swapping the endpoints of an edge is considered a nontrivial automorphism
of that edge). With this convention, define $\fn{o_E}{E(G)}{\N}$ by 
$o_E(e)=|O(e)|$.
\end{definition}

\begin{definition}
A tree is called {\em perfect of Type $n$} if 
\begin{enumerate}
\item $n=1$: the tree has a single vertex.
\item $n=2$: the tree is a binary tree.
\item $n>2$: the tree consists of $n$ binary trees linked to a common root.
\end{enumerate}
\end{definition}

\begin{lemma}[Product Decomposition]\label{prod}
Suppose $G$ is a tree with an edge $e$ such that $o_E(e)=1$. Removing $e$ 
results in two connected trees $G_1$ and $G_2$. In this case, $\GAut~G$ 
decomposes
as $\GAut~G_1\times \GAut~G_2$. Moreover, if $G$ is optimal, so are $G_1$ and 
$G_2$.
\end{lemma}

\begin{proof}
The first part is clear: $e$ is not moved by any automorphism, so any 
nontrivial automorphism must be an automorphism of $G_1$ or $G_2$ (or a 
composition thereof). The second part is also easy: if $G_1$ is not optimal, 
then $G_1$ (as a subgraph of
$G$) could be replaced by an optimal graph with the same number of leaves of 
$G_1$, contradicting optimality of $G$. 
\end{proof}

The following lemma is the most essential part of the proof. It states that
if a vertex $v$ is moved by the automorphism group, then its branches should
all be isomorphic (except along the geodesic leading to the absolute root). If
not, the various copies of the branches attached to vertices in the orbit of
$v$ should be removed and grouped together to increase symmetry. 

\begin{lemma}[Terminal Symmetry]\label{term}
Suppose $v$ is a vertex with $o_V(v)>1$. Then the branches of $G_v$ around
$v$ are all isomorphic.
\end{lemma}

\begin{proof}
The strategy is this: if $v$ is a moving vertex, and has two non-isomorphic
branches, these branches also move, and a more symmetric graph can be created
by grouping like branches together.

Suppose that there are two vertices $v_1$ and $v_2$ adjacent to $v$ in
$G_v$ such that $G_{v_1}\not\cong G_{v_2}$. Since $v$ moves, there are other
copies of $G_{v_i}$ in the graph $G$. 

Denote by $G_i'$ the tree obtained by removing $G_{v_i}$ and its orbit
under the automorphism group of $G$. Then $\GAut~G$ surjects onto $\GAut~G_i'$
with kernel those automorphisms fixing the vertices of $G$ not in the image of
$G_{v_i}$ under the action of $\GAut~G$. The order of this kernel is 
$|\GAut~G_{v_i}|^{o_V(v_i)}$.

Let $G_i$ denote the stabilization of $G_i'$. Then $\GAut~G_i'$
is isomorphic to $\GAut~G_i$.

Now construct a graph $G_i^\circ$ from $G$ by removing everything from $G$
except the orbit of $G_{v_i}$ and the corresponding geodesics to the absolute
root
(including the root edge if the absolute root is virtual) and stabilizing the 
result.
Since $G_{v_1}\not\cong G_{v_2}$, $G_1$ and $G_1^\circ$ are nontrivial, and
the sum of their numbers of leaves is the number of leaves of $G$. Join
$G_1$ and $G_1^\circ$ at their roots (making an appropriate construction
when the root is virtual, or when this makes the valence too high at the new
root). The resulting graph has more  automorphisms than $G$, since
the order of $\GAut~G_1^\circ$ is at least $o_V(v_1)|\GAut~G_{v_1}|^{o_V(v_1)}$.
This contradicts the optimality of $G$ and proves the lemma.
\end{proof}

\begin{proposition}\label{mins}
Suppose $G$ is optimal. Then $\min(o_E)\leq 5$ and $\min(o_V)\leq 2$. Moreover,
if $\min(o_E)\geq 3$, then the graph is $\min(o_E)$ isomorphic subtrees 
attached to a common root.
\end{proposition}

\begin{proof}
First suppose $G$ has an invariant vertex. Then clearly, $\min(o_V)=1$ is
attained at this vertex. Since an optimal graph has valence at most five,
the orbit of an edge with this root has at most five elements. This shows that
$\min(o_E)\leq 5$, since these five edges may be permuted at most among
themselves.

Now suppose $G$ has no invariant vertex. Then there must be an invariant
edge. This edge is either carried in an oriented way onto itself, or
the orientation is reversed, so $\min(o_E)\leq 2$. Since the edge is
invariant, its endpoints can at most be taken to each other, so
$\min(o_V)\leq 2$ in this case as well. 

The graph has an absolute root or an absolute virtual root. If there is an
absolute virtual root $e$, this edge has $o_E(e)=1$ or $o_E(e)=2$, 
In case $\min(o_E)\geq 3$, there is an absolute root $v$. Consider the 
isomorphism classes of the branches from $v$. If any class has a single
member, the edge from $v$ to the root of that branch is an invariant edge, 
which contradicts $\min(o_E)\geq 3$. Now if there are at least two isomorphism
classes (each with more than one member), we have a contradiction of optimality
by the proof of Lemma
\ref{valred}. It is clear that one of these edges attains
the minimum of $o_E$ (since the whole graph rotates around $v$), which proves
the last statement of the proposition
\end{proof}

\begin{proposition}[Doubling Lemma]\label{double}
An optimal graph with $2g$ leaves can be obtained from an optimal graph
with $g$ leaves by replacing each leaf with a vertex attached to two leaves.
\end{proposition}

\begin{figure}[ht]
\begin{center}
\begin{picture}(130,45)
\put(5,15){\circle*{3}}
\put(5,15){\line(1,0){10}}
\put(15,15){\circle*{3}}
\put(15,15){\line(1,0){15}}
\put(15,15){\line(0,1){10}}
\put(15,15){\line(0,-1){10}}
\put(15,5){\circle*{3}}
\put(15,25){\circle*{3}}
\put(30,15){\circle*{3}}
\put(30,15){\line(1,0){15}}
\put(30,15){\line(0,1){15}}
\put(30,30){\circle*{3}}
\put(30,30){\line(0,1){10}}
\put(30,30){\line(-1,0){10}}
\put(30,30){\line(1,0){10}}
\put(20,30){\circle*{3}}
\put(40,30){\circle*{3}}
\put(30,40){\circle*{3}}
\put(45,15){\circle*{3}}
\put(45,15){\line(1,0){10}}
\put(45,15){\line(0,1){10}}
\put(45,15){\line(0,-1){10}}
\put(55,15){\circle*{3}}
\put(45,5){\circle*{3}}
\put(45,25){\circle*{3}}

\put(60,15){\vector(1,0){10}}

\put(75,10){\circle*{3}}
\put(75,10){\line(1,0){10}}
\put(85,10){\circle*{3}}
\put(85,10){\line(1,0){10}}
\put(85,10){\line(0,1){10}}
\put(95,10){\circle*{3}}
\put(85,20){\circle*{3}}
\put(85,20){\line(0,1){10}}
\put(85,20){\line(1,0){15}}
\put(85,30){\circle*{3}}
\put(85,30){\line(-1,0){10}}
\put(85,30){\line(1,0){10}}
\put(75,30){\circle*{3}}
\put(95,30){\circle*{3}}
\put(100,20){\circle*{3}}
\put(100,20){\line(0,1){15}}
\put(100,20){\line(1,0){15}}
\put(100,35){\circle*{3}}
\put(115,20){\circle*{3}}
\put(115,20){\line(0,1){10}}
\put(115,20){\line(0,-1){10}}
\put(115,10){\circle*{3}}
\put(115,30){\circle*{3}}
\put(115,10){\line(-1,0){10}}
\put(115,10){\line(1,0){10}}
\put(105,10){\circle*{3}}
\put(125,10){\circle*{3}}
\put(115,30){\line(-1,0){10}}
\put(115,30){\line(1,0){10}}
\put(105,30){\circle*{3}}
\put(125,30){\circle*{3}}
\end{picture}
\end{center}
\caption{Splitting a 3-3-3 configuration in the Doubling Lemma}
\label{dub}
\end{figure}
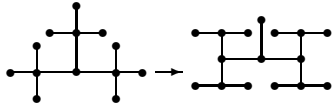

\begin{proof}
It suffices to show that an optimal graph with $2g$ leaves has the property
that exactly two leaves are connected to a vertex which is connected to any
leaves. Such a graph is certainly obtained by ``doubling''. Conversely, if 
such a graph has its leaves removed to obtain a graph with $g$ leaves which is
not optimal, doubling an optimal graph with $g$ leaves will produce a more
symmetric graph with $2g$ leaves.

For reasons of valence, the only configurations of leaves other than two
per branch are branches with one leaf or branches with three, except for
the two cases of four leaves around a root and five leaves around a root.
Five is not even, so the lemma doesn't apply, and there is an optimal tree 
with four leaves which is doubled from the optimal tree with two leaves.
By the Terminal
Symmetry Lemma and stability of the curve in question, branches
with one leaf are not permuted by the geometric automorphism group: by 
stability, there must be at least two edges other than the one to the leaf,
and removing one on the geodesic to the root leaves a tree with at least
one rational branch and one elliptic leaf. 
If there
are an even number of such branches, they may be combined pairwise to increase
the order of the automorphism group (remove one leaf and place it on a branch
with another, yielding an involution). Therefore in an {\em optimal} graph, 
there is at most one such branch.

Now suppose that $v_1$ is a vertex adjacent to three leaves. We claim that
$o_V(v_1)=1$. Denote by $v_0$ the vertex one step from $v_1$ towards the
absolute root (if $v_1$ is the absolute root, the claim is clearly true). 

Suppose that $o_V(v_0)>1$. Then the Terminal Symmetry Lemma implies that
all vertices one unit away from $v_0$ in the tree $G_{v_0}$ are branches with
three leaves. By valence considerations, the only possibilities are that
$G_{v_0}$ has two to five branches. If there are two branches, split the six
leaves into three branches with two leaves each. If there are three, split
into the configuration shown in Figure \ref{dub}. In both cases, the 
contribution to
automorphisms increases, in the first case from 18 to 24, and in the second 
from 81 to 128. Similar constructions can be performed to a configuration of
four or five branches of three leaves around a root (the answers are given
by the Main Theorem) to get more automorphisms. This contradicts optimality, 
so $o_V(v_0)=1$.

Now, if there is another branch of $G_{v_0}$ adjacent to $v_0$ which has three 
leaves, these
could be combined as in the previous paragraph with the leaves around $v_1$,
contradicting optimality. Hence $v_1$ is the only branch of $G_{v_0}$ 
adjacent to $v_0$ with
three leaves, so if it is moved by some automorphism of $G$, $v_0$ will follow.
This contradicts the fact that $o_V(v_0)=1$.

Thus branches with an odd number of leaves do not move around the graph. 
Therefore, they may be broken up to increase symmetry: pairing two branches
with a single leaf adds an involution switching the leaves. Pairing a single
leaf with a branch with three leaves and splitting into pairs increases the
automorphism group by a factor of at least 8/3. Similarly, two branches with
three leaves each can be combined. These constructions contradict optimality
of the graph, and we conclude that an optimal graph of even order has exactly
two leaves on every branch which has any leaves at all. The proposition is
proved.
\end{proof}

\begin{lemma}[Strong Terminal Symmetry]\label{sts}
If $o_V(v)>1$ for some vertex $v$, then the branches of $G_v$ around
$v$ are all isomorphic {\em perfect} trees.
\end{lemma}

\begin{proof}
If the number of leaves of one of these branches is even, then the doubling
lemma allows us to prove the result by induction. If the number of leaves
on a branch is not even, then each branch has a subbranch with three leaves,
but we have seen that such configurations are not optimal, so the number of 
leaves must be even.
\end{proof}

Now the preliminaries are in place, and the Main Theorem may be proved.

\begin{proof}[Proof of Main Theorem]
The genus two case is easily checked by hand. The base case for part two
is that of genus three, which follows from the fact that there is a unique
simple graph among dual graphs of genus three curves. In genus four, there
are two simple dual graphs, both maximally symmetric, one of which satisfies
the form of the theorem. In genus five, it is also easy to find a maximally
symmetric graph among the simple graphs, which is the final case of the Main
Theorem.

The proof proceeds by induction on the number of binary digits of $g$. Suppose
the result is known for $g$ with $m$ or fewer binary digits. The Doubling 
Lemma then shows that if $g$ has $m+1$ binary digits and the last digit is 
zero, the result follows.

So we may suppose that $g$ has $m+1$ digits and the last digit is one. This
implies that $\min(o_E)$ is not two (otherwise there would be an even number
of leaves). If $\min(o_E)\geq 3$, then we are done by Strong Terminal Symmetry
and Proposition \ref{mins}. 

The only remaining possibility is that when $g$ is odd, $\min(o_V)=1$,
that is, there is an invariant edge. There may be several such edges; let $e$
be an invariant edge where the ratio between the number of vertices on the 
large side and small side is maximized. Remove this edge and call the larger
resulting graph $G_1$ and the smaller resulting graph $G_2$. By product
decomposition, $\GAut~G\cong \GAut~G_1\times \GAut~G_2$ and $G_1$ and $G_2$ 
are optimal.

If $G_2$ has only one vertex, then it contributes nothing to $\GAut~G$. 
Therefore, $G$ is obtained from an optimal graph by adding
a vertex. We may add an edge to an existing appendage so that the new 
resulting appendage is maximally symmetric. This will yield the wrong answer
if the appendage grows too large (i.e. becomes a binary tree). But then
the answer has been given by strong terminal symmetry. So the case of $G_2$
a single vertex is done.

Either $G_1$ or $G_2$ has an odd number of leaves. Suppose first that $G_1$
has an odd number.
By the induction hypothesis, $G_2$ is doubled from a graph with
half as many leaves, so it has no vertices of valence (2,1), (3,1), or (1,3).
On the other hand, $G_1$ must have an invariant vertex of one of these three
types. Since the edge connecting this invariant vertex of valence (2,1) or
(3,1) to a leaf must be
invariant, vertices of valence (2,1) and (3,1) do not occur (the ratio
of the number of vertices in $G_1$ to that of $G_2$ was chosen maximal, and
$G_2$ has at least two vertices). Thus $G_1$ has a vertex of valence
$(1,3)$, unique by induction. Considering the subtree of $G$ rooted at
this vertex shows that $G_2$ has at most two leaves, otherwise $G$ could be
divided at the (1,3) vertex to yield a higher weight ratio. In this case, 
however, since the (1,3) vertex of $G_1$ is invariant, this subtree
can be removed, joined with the branch supporting the two leaves of $G_2$,
and the leaves redistributed to increase the order of the automorphism group.

Therefore the smaller graph $G_2$ has an odd number of leaves. Previous 
arguments on $G_1$ show that it is enough to consider the case that $G_2$ has 
three leaves. By induction, we have one of the following
\begin{enumerate}
\item $G_1$ has $3\cdot 2^n+a$ leaves and is of the form given by the theorem.
If $a+3<2^n$, then $G$ fits the form of the theorem: the three leaves of
$G_2$ are part of an appendage. In any case, the appendage of $G_1$ is 
itself a nested collection of maximally symmetric trees of the types given,
so $G_2$ is attached to the last of these: therefore the problem reduces 
to adding a branch with three leaves to an appendage with six leaves: it
is easy to see that any such configuration can be rearranged to give more
automorphisms, so in fact, a $G_2$ with three leaves does not occur in this
case.
\item $G_1$ has $4\cdot 2^n+b$ leaves and is of the form given by the theorem.
If $b+3<2^{n+1}$, then $G$ fits as in part one. An argument similar to that
given above shows that  this border crossing does not
happen in this case either (in this case, adding a branch with three leaves
to one with two does not give an optimal configuration of five leaves.
\item $G_1$ has $5\cdot 2^n$ leaves. It is clear that adding $G_2$ to a 
maximal $G_1$ as given in the theorem will not give a maximally symmetric 
curve (the root may be broken), so this case does not occur.
\end{enumerate}
\end{proof}

\begin{theorem}\label{form}
The order of the automorphism group of a maximally symmetric stable curve of
genus $g$ is
\begin{enumerate}
\item if $g=2$, $72$.
\item if $g=3\cdot 2^n$, $6^g\cdot 2^{g-3}\cdot 6$.
\item if $g=5\cdot 2^n$, $6^g\cdot 2^{g-5}\cdot 10$.
\item $6^g\cdot 2^{N(g)}\cdot\left(\frac38\right)^{k(g)}\cdot\left(\frac12
\right)^{l(g)}$
\end{enumerate}
where $N(g)$, $k(g)$, and $l(g)$ are computed from the binary expansion of 
$g$ as
follows: starting from the left side, look for successive groups of two 
bits starting with one, disregarding any intermediate zeros; $k(g)$ is the 
number of groups of the form $11$, $l(g)$ is the number of groups of the form
$10$, $N(g)=g-1$ if there is a one remaining on the right end after pairing,
and $N(g)=g$ otherwise.
\end{theorem}

\begin{example}
It is worth illustrating the formula with an example. Write 215 in base two
as 11010111. There is a group 11 at the far left, then an intermediate zero,
then a 10, then a 11, and a ``lonely'' 1. So $k(215)=2$, $l(215)=1$, and
$N(215)=214$.
\end{example}

\begin{proof}
Call a number $g$ {\em special} if after the pairing explained in the statement
of the theorem, there is a ``lonely'' 1 left. Clearly an even number is never
special. The optimal graph for a special
$g$ has an isolated leaf; in other odd genera the last pair is 11, so there
is an isolated branch with three leaves.

The following formulas are easily obtained:
\begin{itemize}
\item If $g$ is special, $2g$ and $2g+1$ are not special. Thus $N(g)=g-1$,
$N(2g)=2g$, and $N(2g+1)=2g+1$; $k(2g)=k(g)$, $k(2g+1)=k(g)+1$ and 
$l(2g)=l(g)+1$, $l(2g+1)=l(g)$.
\item If $g$ is not special, $2g$ is not special, but $2g+1$ is. Thus 
$N(g)=g$, $N(2g)=2g=N(2g+1)$, $k(2g)=k(2g+1)=k(g)$, and $l(2g)=l(2g+1)=l(g)$.
\end{itemize}

The proof naturally proceeds by induction on the number of binary digits of
$g$. The base cases $g=2$, $3$, and $5$ are easily checked by hand. 
Suppose the result is known for $n-1$ binary digits
and that $g$ has $n$ binary digits. If the last bit of $g$ is zero, then
the observations above and the Doubling Lemma prove the result.

If $2g+1$ is not special, then $g$ is special. Doubling $G$ produces an
isolated branch with two leaves (i.e. a branch vertex $v$ with $o_V(v)=1$). By
the proof of the Main Theorem, we may go from such a genus $2g$ curve to
a maximally symmetric genus $2g+1$ curve by adding a leaf to the isolated
branch. In light of the formulas above, it is easy to check that this gives
the desired order of the group.

In the case that $2g+1$ is special, $g$ is not special, and the extra leaf
added in passing from $2g$ to $2g+1$ is attached to an invariant vertex, and
hence adds no automorphisms to the tree. Therefore the formula is also true
in this case.
\end{proof}

\section{Description of all maximally symmetric curves}

\begin{figure}[ht]
\begin{center}
$\begin{array}{c@{\hspace{0.5in}}c}
\begin{picture}(40,35)
\put(5,5){\circle*{3}}
\put(5,5){\line(0,1){10}}
\put(5,15){\circle*{3}}
\put(5,15){\line(0,1){10}}
\put(5,15){\line(1,0){15}}
\put(5,25){\circle*{3}}
\put(20,15){\circle*{3}}
\put(20,15){\line(1,0){15}}
\put(20,15){\line(0,1){15}}
\put(20,15){\line(0,-1){10}}
\put(20,5){\circle*{3}}
\put(20,30){\circle*{3}}
\put(20,30){\line(1,0){10}}
\put(20,30){\line(-1,0){10}}
\put(10,30){\circle*{3}}
\put(30,30){\circle*{3}}
\put(35,15){\circle*{3}}
\put(35,15){\line(0,1){10}}
\put(35,15){\line(0,-1){10}}
\put(35,5){\circle*{3}}
\put(35,25){\circle*{3}}
\end{picture}
&
\begin{picture}(40,45)
\put(5,5){\circle*{3}}
\put(5,5){\line(0,1){10}}
\put(5,15){\circle*{3}}
\put(5,15){\line(0,1){10}}
\put(5,15){\line(1,0){15}}
\put(5,25){\circle*{3}}
\put(20,15){\circle*{3}}
\put(20,15){\line(1,0){15}}
\put(20,15){\line(0,1){15}}
\put(20,30){\circle*{3}}
\put(20,30){\line(1,0){10}}
\put(20,30){\line(-1,0){10}}
\put(20,30){\line(0,1){10}}
\put(20,40){\circle*{3}}
\put(10,30){\circle*{3}}
\put(30,30){\circle*{3}}
\put(35,15){\circle*{3}}
\put(35,15){\line(0,1){10}}
\put(35,15){\line(0,-1){10}}
\put(35,5){\circle*{3}}
\put(35,25){\circle*{3}}
\end{picture}
\end{array}$
\end{center}
\caption{Nonuniqueness in genus seven}
\label{seven}
\end{figure}
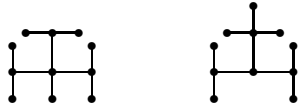

In the previous two sections, we have given one way of finding a maximally
symmetric stable curve of genus $g$. It is natural to ask if the curve found
is unique. The first counterexample occurs in genus four, where there are
two maximally symmetric curves; this happens again in genus seven, pictured 
in Figure \ref{seven}.

As the genus increases, even worse nonuniqueness can occur. However, we can
describe all maximally symmetric curves of a given genus.

\begin{proposition}
A maximally symmetric curve of genus $g$ has
\begin{enumerate}
\item all components $\Pro^1$ or $E$,
\item dual graph a tree with all leaves elliptic and other vertices rational.
\end{enumerate}
\end{proposition}

\begin{proof}
In all of the reduction steps, we actually gain automorphisms with one
exception: reducing valence at a vertex where neighboring vertices are
non-isomorphic does not necessarily add any automorphisms. Therefore, 
conditions of valence are dropped from the conditions of simplicity giving
the present result.
\end{proof}

An inspection of the proof of the Terminal Symmetry lemma shows that the
restrictions on valence were not used there. Therefore the Terminal Symmetry
lemma may be applied to non-simple curves.

\begin{lemma}[Perfect Structure]
Let $G_0$ indicate the subtree of $G$ fixed by all geometric automorphisms of 
$G$. If
$G$ is optimal, then each leaf of $G_0$ is the root of a perfect subtree of
$G$. Moreover, trees of types 4 and 5 may occur only when $G_0$ is a point.
\end{lemma}

\begin{proof}
Suppose $t_0$ is a leaf of $G_0$ which is not the root of a perfect subtree.
Since $t_0$ is a leaf of the fixed subtree, none of its neighbors outside of
the fixed subtree are fixed. Then Strong Terminal Symmetry says that the 
subtrees
whose roots are these neighbors (call then $v_i$) have all isomorphic
branches, and these branches are perfect.

We now claim that the $G_{v_i}$ are all isomorphic, and that there are at
most five of them. The second claim follows from the first, since if all the
branches are isomorphic and greater in number than five, they can be 
rearranged (Lemma \ref{valred}) to contradict the optimality of $G$.

Supposing at least two of the subtrees are non-isomorphic, there are at least
two orbits of trees around $t_0$. However, since the point of attachment of 
$t_0$ to the rest of the invariant tree must be fixed by any automorphism
of the $\Pro^1$ corresponding to $t_0$ in the curve, there are
not enough automorphisms of the $\Pro^1$ left to realize every possible 
graph symmetry (since all orbits must be nontrivial). All of the symmetries can be realized by splitting $t_0$
and rearranging the branches, contradicting optimality. Therefore all of the 
subtrees are isomorphic and there are at most five.

The lemma follows, since these subtrees themselves are perfect.

The last part of the lemma follows since a tree of Type 4 or 5 attached to a 
leaf of a non-trivial $G_0$ cannot realize its full symmetry group . Then 
splitting the tree in two trees increases the order of the automorphism group, 
contradicting the optimality of $G$. 
\end{proof}



The following definition, especially the last condition, serves to isolate
the binary pairs 10 and 11 occurring in the proof of the Main Theorem. The
exception in the last item is necessary to note: without it, there is
no strict optimal graph in genus eleven. The upshot of the definition here
and the proofs below is that the behavior in genus seven and eleven somehow
is the whole picture.

\begin{definition}
Call a tree $G$ {\em strict optimal} if
\begin{enumerate}
\item $G$ is the union of the fixed subtree $G_0$ and $k$ perfect subtrees
$G_i$ whose roots are on $G_0$; index the $G_i$ so that their numbers of
leaves $N_i$ are decreasing;
\item the roots of the $G_i$ are the leaves of $G_0$;
\item the valence at interior vertices of $G_0$ is exactly three 
\item $N_i\geq 4N_{i-1}$ for $i=2,\ldots,k$ except when $G_i$ is Type 2
with $2^{s+3}$ and $G_{i-1}$ is Type 3 with $3\cdot 2^s$ leaves.
\end{enumerate}
\end{definition}

\begin{lemma}
A strict optimal tree is optimal.
\end{lemma}

\begin{proof}
This follows from Product Decomposition (the automorphism group of a strict
optimal tree is the direct product of the automorphism groups of its subtrees
rooted at leaves on the invariant tree) and the last condition in the
definition of strict optimal, which allows us to compute the order of the
automorphism group of a strict optimal tree and see that it has the value
given in \ref{form}.
\end{proof}

\begin{definition}
If an optimal graph has two perfect subtrees $G_i$ and $G_j$ 
with $2^{s+2}$ and $3\cdot 2^s$ leaves, respectively, we define a {\em neutral
move of Type I}: remove a perfect subtree with $2^s$ leaves from $G_j$ (leaving it a
binary tree with $2^{s+1}$ leaves) and attach it to a different vertex of 
$G_0$, splitting an edge with a new vertex if necessary to keep valence low
(note in particular that a neutral move for a given tree is not unique). 
Now attach the rest
of $G_j$ (the aforementioned binary tree) at the root of $G_i$, resulting in
a perfect subtree with $3\cdot 2^{s+1}$ leaves. This process will tacitly
be followed by any stabilization of the graph necessitated by bad choices.

We define a {\em neutral move of Type II} only in case the (optimal) tree has
precisely $2^n$ ($n\geq 2$) leaves. In that case, if four perfect binary 
trees are sharing the common root $G_0$, we separate two of the four trees
around another rational node, linked by an edge to the original rational root. 
\end{definition}

This definition is easier to grasp with the examples of nonuniqueness in
genus seven in mind. The left hand example in Figure \ref{seven} is a strict
optimal tree. Its invariant subtree is the lower central segment, bearing
a perfect subtree with six leaves, and one with a single leaf. This satisfies
the inequalities in the definition of strictness. The right hand example has
its most central vertical segment as invariant subtree, bearing perfect
subtrees with four and three leaves. This violates strictness. Here we
are in the situation of the previous definition with $s=0$. Remove the 
highest vertical edge in the figure and place it at the lower vertex of
the invariant tree. All neutral moves are obtained by ``doubling''
this move.

\begin{proposition}
A neutral move preserves the order of the geometric automorphism group of
the graph.
\end{proposition}

\begin{proof}
Using the Doubling Lemma backwards, the situation of the definition of 
neutral move reduces to the case of subtrees of orders three and four, 
where the explanation of the genus seven example clearly shows that the
automorphism group does not grow or shrink.
\end{proof}

The following theorem shows that the Main Theorem is close to giving all
maximally symmetric curves. The motivation is trying to reverse the formula
of Theorem \ref{form}. Pairing binary digits, we try to reconstruct the
tree. A neutral move occurs when a odd-length sequence of ones (three or more)
occurs in the binary expansion. The proof follows slightly different lines.

\begin{theorem}
Every maximally symmetric genus $g$ curve has either a strict optimal dual
graph, or its dual graph can be made strict optimal by a sequence of neutral
moves, valence reduction, and stabilization. 
\end{theorem}

\begin{proof}
Clearly a maximally symmetric curve must have an optimal graph. As previously,
denote the subtrees rooted at leaves of the invariant tree $G_0$ by
$G_1,\ldots,G_k$, ordered so that the number of leaves in these subtrees is 
decreasing. Since the $G_i$ are rooted at invariant nodes, no two have the
same number of leaves: if they did, since they are perfect, they would be
isomorphic, and the graph could be rearranged to be more symmetric. Therefore
$N_1>N_2>\cdots>N_k$. Also, we have $\GAut~G=\prod_{i=1}^{k}\GAut~G_i$.

If $G$ itself is a perfect tree, it can have $2^n$, $3\cdot 2^n$ or 
$5\cdot 2^n$ leaves. In the last two cases, or in the first when $n=1$, 
there are no other optimal graphs possible (by Strong Terminal Symmetry). In 
the first case with $n\geq 2$, 
there are two possibilities: either $\min(o_E)=4$ or $\min(o_E)=2$. If the 
first possibility occurs, we do a neutral move of Type II to get to a
strictly optimal graph, while in the second case the graph is already strictly 
optimal.

If $G$ is not a perfect tree, we have at least two distinct leaves in $G_0$. 
Moreover, we know from valence reasons that none of the $G_i$'s can be perfect
of Type 4 or 5.

By induction, we may remove the subtree $G_1$ and assume that the tree 
remaining is optimal and satisfies the conclusion of the theorem.

If $G_1$ is perfect of Type 3, $N_1=3\cdot 2^s$ and either
\begin{enumerate}
\item $G_2$ is of Type 2, $N_2=2^p$. Then $p\leq s+1$; if $p=s+1$ or $p=s$,
we are not optimal (``undouble'' down to the case three versus two or three
versus one and note that we may rearrange). Therefore $N_1\geq 6N_2>4N_2$.
\item $G_2$ is of Type 3, $N_2=3\cdot 2^p$. Then $p\leq s-1$, and an 
undoubling argument shows that if $p=s-1$, then $G$ is not optimal as in
the first case. Therefore again, $N_1\geq 4N_2$.
\end{enumerate}
Thus if $G_1$ is of Type 3, $G$ is strict optimal by induction.

Now if $G_1$ is Type 2 with $N_1=2^{s+2}$ leaves (the case $g=2$ is clear),
again, either
\begin{enumerate}
\item $G_2$ is of Type 2; $N_2=2^p$ and $p<s+1$ (by optimality), so
$N_1\geq 4N_2$.
\item $G_2$ is of Type 3; $N_2=3\cdot 2^p$. This divides into further cases:
\begin{enumerate}
\item $p<s-1$: $N_1\geq 4N_2$.
\item $p=s-1$: $N_1= \frac83N_2$: this is the ``exception'' in the condition
of strictness.
\item $p=s$: do a neutral move to change $N_1$ to $3\cdot 2^{s+1}$ and
$N_2$ to $2^s$. Then in the new tree, $N_1\geq 4N_2$.
\end{enumerate}
\end{enumerate}

Having achieved the numerical condition, it is easy to achieve valence three
at every interior vertex of the fixed tree. If a perfect subtree is attached
to an interior node, it may be branched out so it is rooted at a leaf (the
new edge will also be invariant).
\end{proof}

\begin{remark}
In many cases, the number of maximally symmetric curves is finite (in some
cases, notably $5\cdot 2^n$, $3\cdot 2^n$ and $2^n$ it is unique). 
But there are cases
where there is a positive dimensional family of maximally symmetric curves
(exactly when $k(g)+l(g)+g-N(g)-3$ is positive, in which case this quantity
is the dimension of the family of maximally symmetric stable curves).
The easiest example to see is probably in genus $1+4+16+64$. By the Main
Theorem, a maximally symmetric curve is constructed by first arranging
a binary tree with 64 leaves, then attaching to its root a maximally symmetric
genus $1+4+16$ curve, and so on. In the end, there is a root connected to
binary trees with one, four, sixteen, and sixty-four leaves. This root
could be split, as in the previous theorem to yield strict optimal trees
(note in particular that strict optimal trees are not unique in a given genus),
but this does not affect automorphisms, so we might as well keep all four
branches tied to a single root. However, the automorphism group of
$\Pro^1$ is only three-point transitive, so after attaching the first three
branches, there are infinitely many choices for the point of attachment of
the fourth branch.
\end{remark}

\end{document}